\documentclass[review]{elsarticle}
\usepackage{lineno,hyperref}
\usepackage{amsmath,amssymb,amsfonts}
\usepackage[ruled,vlined,linesnumbered]{algorithm2e}
\usepackage{graphicx}
\usepackage{textcomp}
\usepackage{xcolor}
\usepackage{multirow}
\usepackage{threeparttable}
\usepackage{subfigure} 
\usepackage{booktabs} 
\usepackage{amsmath}
\usepackage{multirow}  
\usepackage{rotating} 
\usepackage{mathrsfs} 
\usepackage{amssymb}  
\usepackage{natbib}
\usepackage{nomencl}
\usepackage{etoolbox}
\usepackage{tabularx} 

\usepackage{tikz}

\usepackage{tabularx} 
\usepackage{array}    

\definecolor{c1}{rgb}{0,0,1}

\usepackage{geometry}
\makeatletter
\def\ps@pprintTitle{%
	\let\@oddhead\@empty
	\let\@evenhead\@empty
	\let\@oddfoot\@empty
	\let\@evenfoot\@oddfoot
}
\makeatother

\usepackage{numcompress}

\begin{document}

\begin{frontmatter}

\title{Battery swapping station location for electric vehicles: a simulation optimization approach}

\author{Guangyuan Liu}
\address{School of Software Engineering, Tongji University, Shanghai, China}


\author{Yu Zhang, Tianshi Ming}
\address{School of Electronic and Information Engineering, Tongji University, Shanghai, China}

\author{Chunlong Yu}
\address{School of Mechanical Engineering, Tongji University, Shanghai, China}
\cortext[mycorrespondingauthor]{Corresponding author}
\ead{chunlong_yu@tongji.edu.cn}





\begin{abstract}

Electric vehicles (EVs) face significant energy supply challenges due to long charging times and congestion at charging stations. Battery swapping stations (BSSs) offer a faster alternative for energy replenishment, but their deployment costs are considerably higher than those of charging stations. As a result, selecting optimal locations for BSSs is crucial to improve their accessibility and utilization. Most existing studies model the BSS location problem using deterministic and static approaches, often overlooking the impact of stochastic and dynamic factors on solution quality. This paper addresses the facility location problem for BSSs within a city network, considering stochastic battery swapping demand. The objective is to optimize the placement of a given set of BSSs to minimize demand loss. To achieve this, we first develop a mathematical programming model for the problem. Then, we propose a simulation optimization method based on a large neighborhood search framework to handle large-scale instances. To reduce the computational cost of simulations, Bayesian optimization is employed to solve the single-station allocation subproblem during the repair process. Numerical experiments demonstrate the efficiency of the proposed approach and highlight the importance of incorporating dynamic factors in decision-making.
\end{abstract}


\end{frontmatter}

\section{Introduction}

EVs are increasingly viewed as a key part of the future of urban transportation. Compared to traditional internal combustion engine vehicles, EVs are powered by renewable electricity, resulting in lower environmental pollution. They also offer higher energy efficiency and reduced operating costs. Additionally, their simpler drive systems and greater potential for automation make EVs well-suited for autonomous driving technologies.

Despite their advantages, one of the main challenges to wider EV adoption remains energy supply, primarily managed through charging and battery swapping. Charging can be time-consuming, often taking 30 minutes or more, leading to potential congestion at public EV charging stations. While charging addresses routine energy needs, it may not fully meet the demands of long-distance travel or emergencies, which can contribute to range anxiety and reduce customer satisfaction.

BSSs can replenish energy within a few minutes, greatly reducing customer waiting time and offering promising potential for future development. Moreover, the widespread adoption of swapping stations could ease the financial burden of costly batteries, as customers would be able to rent batteries for a low service fee rather than purchase them outright. Swapping stations also help reduce strain on the electrical grid, as batteries can be charged off-site. In addition, battery swapping allows for the integration of the latest battery technology, encouraging the use of more cost-effective and environmentally friendly options.

However, the construction and operational costs of BSSs are significantly higher than those of EV charging stations (EVCSs), requiring substantial investments in parking spaces, mechanical equipment, battery storage, and personnel. The number of BSSs is generally limited in a city compared to the EVCSs. Therefore, identifying suitable locations for swapping stations can enhance accessibility to battery swapping services for consumers, while potentially reducing construction and configuration costs for BSS operators.



Existing research on EVCS location primarily uses mathematical programming to build simplified models of real systems, optimizing them through linear programming methods. For example, Lu Fang's research categorized EVCSs into two types: those on highways and those in urban areas, each corresponding to path-based and point-based charging needs, respectively. However, these studies often overlook the overall optimality of the complex traffic network and base the average battery swapping time on queueing theory, which lacks the realism of simulation systems discussed later.

The location problem for swapping stations differs markedly from traditional fuel stations or charging piles due to their operational characteristics like short service times and multiple failure modes. Additionally, the existing prevalence of charging piles means that EV users will opt for charging when swapping needs are unmet, a previously nonexistent issue. Moreover, user demand can be categorized into point and path demands, with swapping stations better addressing the latter, which has different location decision implications.

Recent research on EV BSS location is still limited. Studies have focused on location planning, considering user behavior capabilities, with the aim of minimizing total costs and maximizing user satisfaction and energy efficiency during battery swaps. These studies employ multi-objective management planning models and solve them using methods like YALMIP/CPLEX, analyzing the relationship between key parameters and outcomes. However, these mathematical programming models often face non-convex constraints that require transformation for solution feasibility.

Our extensive literature review reveals that most swapping station location research simplifies the problem using mathematical programming methods, which fails to capture the complexity of real-world constraints. For example, users consider multiple factors when choosing a swapping station, such as detour distance, proximity to their current location, the maximum distance their vehicle can travel with the remaining charge, and typical congestion at the station. Additionally, urban road networks are vast, often comprising thousands of intersections and connecting roads, with numerous high-frequency path demands each following a different probability distribution. Furthermore, the mechanical complexity of swapping stations inevitably leads to failures, disrupting operations until repairs are completed.

To address these complexities, we propose a novel approach that combines simulation models with optimization algorithms (simulation optimization). This approach allows for the realistic representation of complex systems without the oversimplifications common in mathematical programming models, eliminating non-convex constraint issues and avoiding the need for queueing theory models. This approach can handle large-scale problems and provides real-time monitoring and feedback.

This paper is structured as follows: Chapter 2 reviews relevant literature, Chapter 3 describes the problem, Chapter 4 introduces our methodology, Chapter 5 presents numerical experiments, and Chapter 6 concludes the paper.

\section{Literature review}

At present, many scholars have conducted relevant research on the site selection of charging and battery swapping facilities based on the classic site selection model and the characteristics of charging and battery swapping facilities.

\subsection{EVCS location problem}
The location of EVCS is a crucial aspect in the development and implementation of EVs. Various studies have been conducted to address the challenges and complexities involved in determining the optimal locations for EVCS. \citet{liu2018allocation} focused on factors such as charging satisfaction and distributed renewables integration in the allocation optimization of EVCS. They constructed a multi-objective function considering voltage fluctuation, load fluctuation, and connected capacity of energy storage in EVCS. \citet{xu2018interval} introduced the entropy weight method and a rough consensus reaching process to ensure rational decision-making in EVCS location, proposing a multi-criteria group decision-making framework with linguistic information. \citet{ju2019study} proposed an effective comprehensive framework for evaluating and selecting optimal EVCS locations under a picture fuzzy environment. \citet{hosseini2019development} introduced a Bayesian Network model for optimal site selection of EVCS, considering both quantitative and qualitative factors, providing a new research perspective on uncertainty and qualitative factors in decision-making. \citet{liu2019coordinated} developed an integrated multi-criteria decision-making approach using grey DEMATEL and UL-MULTIMOORA for determining the most suitable EVCS site based on multiple criteria. \citet{zhang2019locating} considered service risk factors in EVCS site selection, introducing the Improved Whale Optimization Algorithm for the model. \citet{wei2020algorithms} investigated the grey relational analysis method for probabilistic uncertain linguistic MAGDM in EVCS site selection, emphasizing the importance of site selection in the whole life cycle of EVCS construction. \citet{feng2021novel} proposed a novel multi-criteria decision-making method for selecting EVCS sites from a sustainable perspective, highlighting the stability and reliability of the method. \citet{ghosh2021application} applied hexagonal fuzzy multi-criteria decision-making 
(MCDM) methodology for EVCS site selection, incorporating geographic information systems with MCDM techniques. \citet{mishra2021single} introduced a single-valued neutrosophic similarity measure-based additive ratio assessment framework for optimal site selection of EVCS, showcasing the methodology's performance in evaluating EVCS sites under a single-valued neutrosophic environment. Overall, these studies demonstrate the importance of considering various factors, such as charging satisfaction, renewables integration, service risk, and qualitative aspects, in the site selection EVCSs to ensure efficient and sustainable deployment.

The literature on EVCS site selection encompasses various methodologies and approaches to optimize the location of these stations. \citet{meng2011optimization} introduced game theory as a means to optimize the location of EVCS, building an optimization game model for this purpose. \citet{tang2013optimal} utilized a fuzzy analysis and Analytical Hierarchy Process (AHP) method to evaluate comprehensive benefits in site selection, incorporating expert consultation. \citet{meng2013evaluation} proposed a AHP for evaluating EVCS sitting programs.  \citet{tao2013location} suggested a site selection method based on the Floyd shortest path algorithm to minimize total charging distance in a region. \citet{guo2015optimal} employed a MCDM method from a sustainability perspective, considering environmental, economic, and social criteria for site selection. \citet{wu2017decision} developed a decision framework for EVCS site selection in residential communities using triangular intuitionistic fuzzy numbers (TIFNs). \citet{chen2017optimal} focused on wind-solar complementary power generation projects for large-scale EVCSs to reduce the impact on the power system.  \citet{al2020case} investigated optimal energy management for a fast-charging station in Qatar with multiple renewable energy and storage systems. Overall, these studies highlight the importance of considering various factors and methodologies in selecting optimal sites for EVCSs.


\subsection{BSS location problem}

The site selection of  EV BSSs is a critical aspect of infrastructure planning to support the growing adoption of EVs. Several studies have focused on developing frameworks and methodologies to aid in the decision-making process for site selection. \citet{wang2020study} proposed a site selection framework for BSS based on a MCDM method, providing investors with a critical tool to choose the most suitable alternative. \citet{WANG2020102149} proposes a site selection framework considering three criteria: land occupation cost, driver comfort, and impact on power grid load levels. A fuzzy decisionmaking trial and evaluation laboratory method is designed to determine the weights of the three criteria. Then, a fuzzy-based method is adopted to rank the candidate locations. A case in Beijing, China, illustrates the effectiveness and robustness of the proposed BSS location decision framework. Another crucial factor for BSS construction is the routing from EVs to stations. \citet{HOF2017102} introduces a BSS location-routing problem for selecting a BSS location from candidate sites and minimizing the sum of construction and routing costs. I, a \cite{su8060557}, a bi-level optimization framework was developed, in which the upper layer
was responsible for deriving the contract price to the lower layer. The lower layer was to minimize the operation costs of distribution
companies and to return the dispatched power of BSSs to
the upper layer. The optimal solution was obtained using the Mesh
Adaptive Direct Search method. To simplify the problem, it determined
the location and sizing of BSSs separately as their optimality
might be impossible to achieve at the same time. But the efficacy of
the proposed method hasn't been validated in real-world scenarios.
A Cost-Concerned Model and a Goal-Driven Model for highway networks
were respectively developed to determine the optimal location of
BSSs on a single path \cite{MakHo-Yin}.

Due to the development of intelligent EVs and data science,
a data-driven location selection model for BSSs is proposed
in \cite{Data-Driven}. In this paper, GPS location data and electricity requests
are collected for a metropolitan area. Then, a location selection
model is proposed with the following three steps: a hidden
Markov model for map matching and trajectory extraction,
an electricity consumption rate model for demand estimation,
and a clustering strategy for location determination. A realistic
case study involving 13,700 taxis in the area of Shanghai,
China, is presented, and the results outperform those of the
state-of-the-art baseline models.

\subsection{Research gap}

In summary, first, existing mainstream research primarily focuses on the issue of EV charging stations, while there is a lack of studies on the location of battery swap stations. Given the significant differences in operational characteristics and real-world contexts between charging stations and battery swap stations, it is crucial to investigate the location selection for battery swap stations. Second, most existing research is based on queuing theory to derive certain parameter values, which inevitably leads to discrepancies with reality. Moreover, many studies simplify various constraints, making them less applicable to real-world scenarios. Third, nearly all studies on battery swap station location use mathematical programming methods, and the complexity of the models naturally leads to the establishment of multiple constraints, some of which may be non-convex. These constraints must first be converted, and to improve the efficiency of solving such complex constraints, accelerated algorithms must also be explored.

Therefore, it is necessary to extend the research on charging stations to the study of battery swap stations, particularly by leveraging more effective modeling approaches that cater to the specific characteristics of battery swap stations. This will enable a systematic simulation of location selection, utilizing simulation capabilities to test existing data and effectively optimize the models.

\section{Problem formulation}  \label{problem description}
\label{section:prob_description}
In this section, we first introduce the problem characteristics. Then, we formulate the objective function based on queueing theory. Finally, a mix-integer mathematical programming (MILP) model of the problem is given.

\subsection{Problem description}  
A city is represented by a connected graph $G = \{V, E\}$, where $V$ is the set of nodes and $E$ is the set of edges. EVs travel within the city, and each EV that requires battery swapping generates a path demand. The path demand of an EV is represented by its origin-destination (OD) pair $(i_s, i_e)$, where $i_s$ is the origin node and $i_e$ is the destination node. We assume that each EV with a path demand $(i_s, i_e)$ must detour to a BSS within the city to recharge before continuing to its destination. The maximum allowable detour distance is $d_{max}$ for all EVs.

Each BSS, denoted by node \( j \), is modeled as a single-server queue with a mean service rate of \( u_j \) and a finite queue capacity \( K \). While the battery swapping time is typically constant, mechanical failures may occasionally occur, causing the service time to become a random variable following a specified distribution.

We consider two situations in which an EV, say \( i \), may abandon the battery swapping service and opt for charging instead. The first, referred to as Type I path demand loss, occurs when the EV arrives at a BSS, but the queue is full. The second, called Type II path demand loss, arises when no BSS is within the reachable range of EV \( i \) given its battery level \( b_i \), or when the detour distance to any BSS exceeds the maximum allowable detour distance \( d_{max} \).

A BSS can be deployed at any node \( j \in V_a \), where \( V_a \subseteq V \) represents the subset of nodes where construction is feasible. For simplicity, we assume the construction costs are equal for each node in \( V_a \). The objective of the BSS location problem is to select \( N \) nodes from \( V_a \) to deploy the BSSs, in order to minimize the total path demand loss over the long run.

Other features and assumptions of the problem are as follows:
\begin{enumerate}
	\item Each EV, say $i$, has an initial battery level $b_i \in BL$, where $BL$ is a finite set of battery levels within $[0, 1]$. 
	\item All EVs have the same energy consumption per unit distance, denoted by $c$.
	\item A path demand $(i_s, i_e)$ occurs according to a possion process with an arrival rate $\lambda_{(i_s, i_e)}$.
	\item An EV $i$ would select a BSS to perform the battery swapping service according to his/her preference, which relates to the detour distance to reach a BSS, as well as the empirical queue waiting time of a BSS. 
\end{enumerate}


\subsection{Formulation of the objective function}
Let $y_j$ be the binary variable indicating whether a BSS is located at node $j$, a BSS location layout can be represented as $\mathscr{L} = \{y_j, j \in V\}$. To optimize $\mathscr{L}$, we build a closed-form expression of the total demand loss based on queueing theory. 

First, given the arrival rate of each path demand, we derive the customer arrival rate of each BSS in the city network. Let $dist(i_s, i_e)$ be the shortest distance between the origin node $i_s$ and the destination node $i_e$ of a path demand, and $d_{\text{max}}$ be the maximum detour of all vehicles, i.e., the users will visit a node $j$ for swapping battery only if the total distance after visiting node $j$ and then to the destination $i_e$ is less than $dist(i_s, i_e) + d_{\text{max}}$. Using this relationship, we can obtain the set of admissible station nodes of a path demand $(i_s, i_e)$ as follows:
\[
AN_{(i_s,i_e)} = \left\{ j \in J \mid dist(i_s, j) + dist(j, i_e) \leq dist(i_s, i_e) + d_{\text{max}} \right\}.
\]
Similarly, the set of admissible station nodes of a path demand with battery level $b$ is:
\[
AN_{(i_s,i_e)}^b = \left\{ j \in J \mid dist(i_s, j) + dist(j, i_e) \leq dist(i_s, i_e) + d_{\text{max}} \text{ and } dist(i_s, j) \leq \frac{b}{c} \right\},
\]
where $c$ is the energy consumption per distance. 

Let $AI_j$ be the set of path demands for which node $j$ is within their admissible station nodes, i.e.,
$AI_j = \left\{ (i_s, i_e) \mid j \in AN_{(i_s,i_e)} \right\}$. Then, the customer arrival rate of a BSS at node $j$ can be calculated as:
\[
\lambda_j = \sum_{(i_s, i_e) \in AI_j} \lambda_{(i_s,i_e)} \cdot \sum_{b \in BL} p_{(i_s,i_e)}^b \cdot P_{j,(i_s,i_e)}^b ,
\]
where $p_{(i_s,i_e)}^b$ is the probability that a path demand $(i_s, i_e)$ has a battery level $b$, and  $P_{j,(i_s,i_e)}^b$ is the probability that an EV user, say $i$, with initial battery level $b$ having an OD pair of $(i_s, i_e)$ would select the BSS $j \in AN_{(i_s,i_e)}^b$. 

The probability $P_{j,(i_s,i_e)}^b$ is formulated as follows. Let $\epsilon$ be the pre-specified probability guarantee and $l_j$ be the empirical mean queue waiting time of BSS $j$. For a user having an OD pair $(i_s, i_e)$, the utility of swapping battery in BSS $j$ is:
\[
U_{j,(i_s, i_e)} = -\alpha_1 \cdot (dist(i_s, j) + dist(j, i_e) - dist(i_s, i_e)) + \alpha_2 \cdot \epsilon - \alpha_3 \cdot l_j + \epsilon_{qj}. 
\]
Here, the non-negative coefficients $\alpha_i$ ($i=1, 2, 3$) denote the sensitivity to detour distance, pre-specified probability guarantee $\epsilon$, and mean queue waiting time $l_j$. These coefficients can be estimated empirically. The random term $\epsilon_{qj}$ contains all unobserved factors, which are considered independently and identically distributed for different EV users. According to the utility maximization rule, the probability that an EV user $i$ with initial battery level $b$ having an OD pair of $(i_s, i_e)$ would select the BSS $j$ is:
\[
P_{j,(i_s,i_e)}^b = \begin{cases} 
\frac{e^{\pi_{j,(i_s,i_e)}^b \cdot y_j}}{\sum_{k \in AN_{(i_s,i_e)}^b} e^{\pi_{k,(i_s,i_e)}^b} \cdot y_k} & \forall j \in AN_{(i_s,i_e)}^b \\
0 & \text{otherwise}
\end{cases},
\]
where
$\pi_{j,(i_s,i_e)}^b = -\alpha_1 \cdot (dist(i_s, j) + dist(j, i_e) - dist(i_s, i_e)) + \alpha_2 \cdot \epsilon - \alpha_3 \cdot l_j  \text{ }(j \in AN_{(i_s,i_e)}^b).$ \label{pi}

\vspace{1em}

Second, we calculate the demand loss rate. Assume that the service time of each BSS follows a exponential distribution, thus, each BSS can be considered as an M/M/1/K system. For Type I loss, according to the M/M/1/K queueing model, the customer lost rate of BSS $j$ is:
\[
\lambda_j^{\text{lost}} = \lambda_j P_j^{\text{lost}} = \lambda_j \frac{\rho_j^{LL} (1 - \rho_j)}{1 - \rho_j^{LL + 1}}
\]
where $\rho_j = \frac{\lambda_j}{\mu_j}$, and $\mu_j$ is the service rate of BSS $j$. While for Type II loss, the loss rate for the path demand $(i_s, i_e)$ with an initial battery level $b$ is:
\[
\lambda_{(i_s,i_e),b}^{\text{lost}} = \begin{cases}
0 & \text{if } \exists j \in AN_{(i_s,i_e)}^b \text{ such that } y_j = 1, \\
\lambda_{(i_s,i_e)} \cdot p_{(i_s,i_e)}^b & \text{otherwise}.
\end{cases}
\]
Thus, the total Type II loss rate for path demand $(i_s, i_e)$ is
\[
\lambda_{(i_s,i_e)}^{\text{lost}} = \sum_{b \in BL} \lambda_{(i_s,i_e),b}^{\text{lost}}.
\]
In summary, the overall percentage of demand loss is:
\[
\eta_{\text{lost}} = \frac{\sum_{(i_s, i_e) \in D} \lambda_{(i_s,i_e)}^{\text{lost}} + \sum_{j \in J} \lambda_j^{\text{lost}}}{\sum_{(i_s, i_e) \in D} \lambda_{(i_s,i_e)}},
\]
where $D$ is the set of all path demands. 

\subsection{MILP Model}  

To clarify the mathematical model, notations are listed as follows. 

\noindent The sets are:\\
\begin{tabularx}{\textwidth}{ p{2.5cm}X }
	$V:$ & Set of nodes representing the city network;\\
	$J:$ & Set of candidate nodes to locate a BSS, $J \subset V$;\\
	$D:$ & Set of all path demands, where a path demand is represented by an origin-destination pair;\\
	$BL:$ & Set of all possible initial battery levels of the vehicles.\\
	$AN_{(i_s,i_e)}:$ & The set of admissible station nodes of the path demand $(i_s, i_e)$;\\
	$AN_{(i_s,i_e)}^b:$ & The set of admissible station nodes of the path demand with battery level $b$;\\
	$AI_j:$ & The set of path demands for which node $j$ is within their admissible station nodes;\\
\end{tabularx} 

\nomenclature[B]{\(j, i, k\)}{Index for nodes}

\noindent The parameters are:\\
\begin{tabularx}{\textwidth}{ p{2.5cm}X }
	$\lambda_{(i_s,i_e)}:$ & The arrival rate of a path demand originating from node $i_s$ to node $i_e$;\\
	$dist(i_s, i_e):$ & The shortest distance between the starting node $i_s$ and the destination node $i_e$ of a path demand;\\
	$d_{\text{max}}:$ & The maximum allowable detour of all vehicles;\\
	$c:$ & The energy consumption per unit distance;\\
	
	$\epsilon:$ & The pre-specified probability guarantee;\\
	$l_j:$ & The empirical mean waiting time of BSS $j$;\\
	$\pi_{j,(i_s, i_e)}^b:$ & The utility of choosing a BSS $j$ for a user that has a path demand $(i_s, i_e)$ and initial battery level $b$;\\
	$\mu_j:$ & Mean service rate of BSS $j$;\\
	$\mu_j^T:$ & Mean service time of BSS $j$, $\mu_j^T = \frac{1}{\mu_j}$;\\
	$LL:$ & The maximum queue length that EV users are willing to tolerate at a BSS;\\
	$p_{(i_s,i_e)}^b:$ & The probability that a path demand $(i_s, i_e)$ has a battery level $b$;\\
	$N:$ & Total number of BSSs.
\end{tabularx}\\

\noindent The decision variables are:\\
\begin{tabularx}{\textwidth}{ p{2.5cm}X }
	$y_j:$ & 1 if a BSS is located at node $j$, 0 otherwise;\\
\end{tabularx}\\

\noindent The auxiliary variables are:\\
\begin{tabularx}{\textwidth}{ p{2.5cm}X }
	$P_{j,(i_s,i_e)}^b:$ & The probability that a EV user associated with path demand $(i_s, i_e)$ and initial battery level $b$ selects the BSS at node $j$ for battery swapping;\\
	$\lambda_j:$ & The customer arrival rate of the BSS at node $j$;\\
	$\rho_j:$ & Traffic intensity of the BSS located at node $j$, defined as $\rho_j = \frac{\lambda_j}{\mu_j}$;\\
	$\lambda_j^{\text{lost}}:$ & Type I lost rate of a BSS at node $j$;\\
	$\lambda_{(i_s,i_e),b}^{\text{lost}}:$ & Type II lost rate for path demand $(i_s, i_e)$ with an initial battery level $b$;\\
	$\lambda_{(i_s,i_e)}^{\text{lost}}:$ & Total lost rate for path demand $(i_s, i_e)$;\\
	$\eta_{\text{lost}}:$ & Overall percentage of demand loss.
\end{tabularx}\\

The MILP model for the BSS location problem is as follows:
\begin{equation} \label{obj}
    \min \frac{\sum_{(i_s, i_e) \in D} \lambda_{(i_s,i_e)}^{\text{lost}} + \sum_{j \in J} \lambda_j^{\text{lost}} }
{\sum_{(i_s,i_e) \in D} \lambda_{i_s,i_e}} 
\end{equation}

subject to:

\begin{equation} \label{equat:1}
    \lambda_j = \sum_{(i_s, i_e) \in AI_j} \lambda_{(i_s,i_e)} \cdot \sum_{b \in BL} p_{(i_s,i_e)}^b \cdot P_{j,(i_s,i_e)}^b ,\quad \forall j \in J
\end{equation}


\begin{equation} \label{equat:2-1}
	P_{j,(i_s,i_e)}^b = 
		\frac{e^{\pi_{j,(i_s,i_e)}^b} \cdot y_j}{\sum_{k \in AN_{(i_s,i_e)}^b} e^{\pi_{k,(i_s,i_e)}^b} \cdot y_k}, \quad \forall  j \in AN_{(i_s,i_e)}^b, (i_s, i_e) \in D, b \in BL 
\end{equation}

\begin{equation} \label{equat:2-2}
	P_{j,(i_s,i_e)}^b = 0 , \quad \forall j \notin AN^b_{(i_s,i_e)}, (i_s, i_e) \in D, b \in BL
\end{equation}

\begin{equation} \label{equat:3}
    \lambda_{(i_s,i_e),b}^{\text{lost}} \geq \left(1 - \sum_{j \in AN_{(i_s,i_e)}^b} y_j \right) \lambda_{(i_s,i_e)} \cdot p_{(i_s,i_e)}^b ,\quad \forall (i_s, i_e) \in D,  b \in BL
\end{equation}

\begin{equation} \label{equat:4}
    \lambda_{(i_s,i_e)}^{\text{lost}} = \sum_{b \in BL} \lambda_{(i_s,i_e),b}^{\text{lost}} ,\quad \forall (i_s, i_e) \in D
\end{equation}

\begin{equation} \label{equat:5}
    \rho_j = \frac{\lambda_j}{\mu_j} ,\quad \forall j \in J
\end{equation}

\begin{equation} \label{equat:6}
    \lambda_j^{\text{lost}} = \lambda_j \frac{\rho_j^{LL} (1 - \rho_j)}{1 - \rho_j^{LL + 1}} ,\quad \forall j \in J
\end{equation}

\begin{equation} \label{equat:7}
	\sum_{j \in J} y_j = N
\end{equation}

\begin{equation} \label{equat:8}
    y_j \in \{0, 1\} ,\quad \forall j \in J
\end{equation}

Objective function (\ref{obj}) aims to minimize total percentage of demand loss in the city. Constraints (\ref{equat:1}) define the customer arrival rate of each BSS. Constraints (\ref{equat:2-1}) and (\ref{equat:2-2}) give the probability that an EV user would select a specific BSS among all alternatives. Constraints (\ref{equat:3}) and (\ref{equat:4}) calculate the Type II loss rate for each path demand. Constraints (\ref{equat:5}) and (\ref{equat:6}) calculate the Type I loss rate at each BSS. Constraint (\ref{equat:7}) imposes the total number of BSSs established.

\section{Methodology}

In this section, we propose a simulation optimization approach to tackle the problem formulated. 
At each iteration, the simulation model evaluates the demand loss rate of a set of given BSS location plans, and returns the objective function values to the optimization module. While the optimization algorithm iterates and generates the next solutions to evaluate, this procedure terminates until some stopping criteria meets. We build the simulation model in Flexsim software, and the details of model construction is available in \ref{simulation}. In the next subsections, we first introduce the basic framework of the optimization algorithm we adopted, i.e., large neighborhood search, and then propose a neighborhood search approach based on Bayesian optimization to mitigate the computational burden of simulations, and finally the complete algorithm. 




\subsection{Large Neighborhood Search}
\label{Large Neighborhood Search}

Large Neighborhood Search (LNS) \cite{shaw1998using} is a heuristic algorithm commonly used for solving complex combinatorial optimization problems, particularly those with large solution spaces and complex constraints, such as vehicle routing \cite{sacramento2019adaptive}\cite{chen2018adaptive}, scheduling \cite{laborie2007self}\cite{godard2005randomized}, and production planning \cite{praseeratasang2019adaptive}.

\begin{figure}[htb!]
    \centering
    \includegraphics[width=0.63\linewidth]{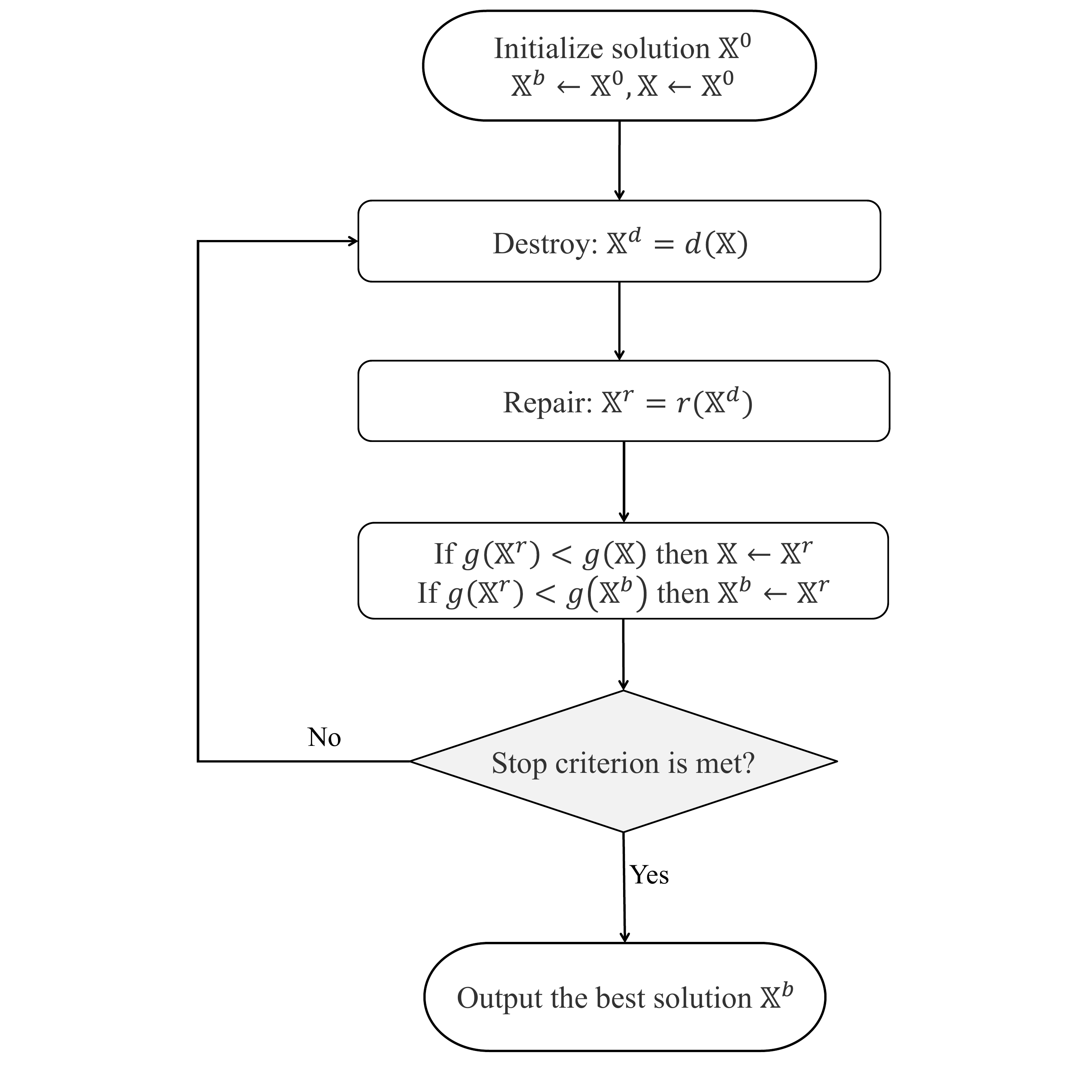}
    \caption{Flowchart of LNS}
    \label{fig:LNS}
\end{figure}

The procedure of LNS is shown in Figure \ref{fig:LNS}. Starting from an initial solution $\mathbb{X}^0$, LNS first applies the destroy method to obtain a partial solution $\mathbb{X}^d$. Then, a repair method is used to repair the partial solution to obtain a complete solution $\mathbb{X}^r$. If the new solution outperforms the current solution, it is accepted and replaces the current solution. Finally, the best solution will be updated if necessary. This procedure iterates until a stopping criterion is met.



In our BSS location problem, we use a simple destroy operator. That is, a set of $k$ BSSs are raondomly selected and removed from the current solution. After, during the repair procedure, these removed BSSs will be re-inserted into the partial solution one by one, until a complete solution is constructed. At each step, the removed BSS will be placed at the best solution such that the objective function is minimized. This indicates the following subproblem:
\begin{equation}
    \min_{\textbf{x} \in \Theta } f( \textbf{x} ), 
    \label{equat:obj}
\end{equation}
where $\textbf{x}$ is the coordinates of the BSS, and 
$\Theta = {\{\textbf{x}_j , \forall j \in J\}} \setminus \mathbb{X}^d $ is the set of possible locations for placing the removed BSS.





To solve (\ref{equat:obj}), the simplest way is to enumerate all possible locations $\textbf{x} \in \Theta$ and select the best position. However, when tackling real-world problems, the simulation model generally involves a large-scale traffic network, thus the simulation would take a considerable amount of time. This makes the enumeration an inefficient optional for solving the subproblem. To deal with the challenge,  we propose a Bayesian Optimization(BO) \cite{frazier2018tutorial} method to solve this subproblem.

\subsection{Bayesian optimization enhanced neighborhood search}



BO is a method for optimizing complex, expensive, or difficult-to-analyze objective functions. 
The basic idea of BO is to construct a metamodel from the observed data points, which guides the search to the promising regions. BO has several advantages: 
(1) High efficiency: by choosing sampling points in a proper way, BO can find promising solutions with relatively few function evaluations. (2) Wide applicability: BO is suitable for non-convex and noisy objective functions, as well as continuous or discrete design spaces. (3) No gradient information: BO does not require gradient information of the objective function and is therefore suitable for black-box functions.

The steps of BO are as folows:  
\begin{itemize}
    \item \textbf{Step 1: Initial design}: Choose a small number of points $\mathbf{X}=\{\textbf{x}_1,\textbf{x}_2,\ldots,\textbf{x}_{n^0} \}$ in the domain, evaluate them to obtain $\mathbf{Y}=\big\{ f(\textbf{x}_1), f(\textbf{x}_2), \ldots, f(\textbf{x}_{n^0}) \big\}$; Set the number of samples $n \leftarrow n^0$. 
    \item \textbf{Step 2: Metamodel construction: } Construct a metamodel, usually, a Gaussian process (GP) model $\mathcal{G} \mathcal{P}(\mu_0 ; k)$ with $\mathbf{X}$ and $\mathbf{Y}$, which predicts the function's behavior across the domain. 

    \item \textbf{Step 3: Sampling:} Employ an acquisition function $\alpha(\cdot)$ to determine the next sample points $\textbf{x}_{n+1} =\arg \max \alpha(\textbf{x})$; set $\textbf{X} \leftarrow  \textbf{X} \cup \{ \textbf{x}_{n+1} \}, \textbf{Y} \leftarrow  \textbf{Y} \cup \{ f(\textbf{x}_{n+1}) \}$;  Update $n \leftarrow n+1$.

    \item \textbf{Step 4: Termination:} If $n$ equals to the maximum number of samples $N$, then stop the algorithm, and output the best sampled solution in $\textbf{X}$, i.e.,
\[ 
\textbf{x}^* = \arg \min_{\textbf{x} \in \textbf{X}} f(\textbf{x})
\]
otherwise go back to \textbf{Step 2}.

\end{itemize}




In the following subsections, we introduce these steps in details. 


\subsubsection{Initial design}
There are many initial strategy to create the initial design, such as random sampling \cite{jones1998efficient}, Latin hypercube sampling \cite{mckay2000comparison} and low-discrepancy sequence sampling \cite{SOBOL196786}. In our case, we use the random sampling to generate $n^0$ points from the domain $\Theta$.




\subsubsection{GP construction}
Suppose we observe $f$ without noise, a GP model $\mathcal{G} \mathcal{P}(\mu_0 ; k)$ is a probabilistic model for the objective function $f$ conditioned on a set of observations $\mathcal{D}_n=\Big\{\big(\textbf{x}_i, f(\textbf{x}_i)\big)\Big\}_{i=1}^n$. It is featured by a mean function $\mu_0$ and a positive-definite covariance function $k$, which is known as the kernel. 

Typically, the mean function is set to 0, and squared exponential kernel is commonly used. More specifically, 
\begin{equation}
    k(\textbf{x}, \textbf{x}') = \sigma_f^2 \exp\left(-\frac{\|\textbf{x} - \textbf{x}'\|^2}{2\ell^2}\right)
\end{equation}
where
    \( \sigma_f^2 \) is the variance parameter of the kernel, controlling the range of the kernel function's output values,
    \( \ell \) is the length scale parameter, which controls the sensitivity of the kernel to the input variation. 
In GP, $\sigma_f^2$ and \( \ell \) is estimated by maximizing the log likelihood 
\begin{equation}
    \log p(\mathbf{f}(\mathbf{X}_n) \mid \mathbf{X})=-\frac{1}{2} \mathbf{f}(\mathbf{X}_n)^\top \mathbf{K}_n^{-1} \mathbf{f}(\mathbf{X}_n) - \frac{1}{2} \log \left| \mathbf{K}_n \right| - \frac{n}{2}\log2\pi,
\end{equation}
where $\mathbf{f}(\mathbf{X}_n)=\left[f\left(\textbf{x}_1\right), \ldots, f\left(\textbf{x}_n\right)\right]^\top$, $\mathbf{X}_n=[\textbf{x}_1, \ldots, \textbf{x}_n]^{\top}$, $\mathbf{K}_n$ is the matrix such that $\left(\mathbf{K}_n\right)_{i j}=k\left(\textbf{x}_i, \textbf{x}_j\right)$.




The resulting prior distribution on $\mathbf{f}(\mathbf{X}_n)$ is
\[
\mathbf{f}\left(\mathbf{X}_n\right) \sim \operatorname{Normal} (\textbf{0}, \mathbf{K}_n ).
\]

Under Gaussian likelihoods, the posterior distribution of $f$ is also a GP,
\begin{align}
    f(\textbf{x}^*) \mid \mathbf{f}(\mathbf{X}_n) {}& \sim \operatorname{Normal}(\mu_n(\textbf{x}^*), \sigma_n^2(\textbf{x}^*))   \label{post distributioin}  \\
    \mu_n(\textbf{x}^*) {}& = \mathbf{k}_n (\textbf{x}^* )^{\top} \mathbf{K}_n^{-1} \mathbf{f}(\mathbf{X}_n)   \label{equation:update bo1}  \\
    \sigma_n^2 (\textbf{x}^* ) {}& =  k (\textbf{x}^*, \textbf{x}^* )-\mathbf{k}_n (\textbf{x}^* )^{\top} \mathbf{K}_n^{-1} \mathbf{k}_n (\textbf{x}^* ),  \label{equation:update bo2}
\end{align}   
where $\mathbf{k}_n (\textbf{x}^* )= [k (\textbf{x}_1, \textbf{x}^* ), \ldots, k (\textbf{x}_n, \textbf{x}^* ) ]^{\top}$ and $\textbf{x}^*$ is the point where the GP is evaluated.

\subsubsection{Acquisition function}
\label{acquisition function}
The acquisition function is used to determine the next point to sample by balancing exploration (sampling uncertain areas) and exploitation (sampling areas likely to yield low function values). To form the acquisition function $\alpha\left(\textbf{x} ; \mathcal{I}_n\right)$, we use the posterior we obtain in \ref{post distributioin}, where $\mathcal{I}_n$ represents the available data set $\mathcal{D}_n$ and the GP structure (kernel, likelihood and parameter values) when $n$ data points are available. The next evaluation is placed at the (numerically estimated) global maximum $\textbf{x}_{n+1}$ of this acquisition function. Many acquisition functions are available now, the most common one is expected improvement.

Let $f_n^*=\min_{m \le n} f(x_m)$ be the optimal choice, where $n$ is the number of points we have evaluated $f$ thus far. Then $f_n^*$ is the previously evaluated point with the minimal observed value.
Now suppose we are evaluating at $\textbf{x}^*$, we will observe $f(\textbf{x}^*)$. After this new evaluation, the improvement in the value of the best observed point is then $\left[f_n^* - f(\textbf{x}^*)\right]^{+}$, where $a^{+}=\max (a, 0)$ indicates the positive part.
While we would like to choose $\textbf{x}^*$ so that this improvement is large, $f(\textbf{x}^*)$ is unknown until after the evaluation. What we can do, however, is to take the expected value of this improvement and choose $\textbf{x}^*$ to maximize it. We define the expected improvement as,
\begin{equation}
    \mathrm{EI}_n(\textbf{x}^*):=E_n\left[\left[f_n^* - f(\textbf{x}^*)\right]^{+}\right],
\end{equation}
where $E_n[\cdot]=E\left[\cdot \mid \mathbf{X}_n, \mathbf{f}(\mathbf{X}_n)\right]$ indicates the expectation taken under the posterior distribution given evaluations of $f$ at $x_1, \ldots x_n$. This posterior distribution is given by \ref{post distributioin}: $f(\textbf{x}^*)$ given $\mathbf{X}_n, \mathbf{f}(\mathbf{X}_n)$ is normally distributed with mean $\mu_n(\textbf{x}^*)$ and variance $\sigma_n^2(\textbf{x}^*)$.

The expected improvement can be evaluated in closed form using integration by parts, as described in \cite{jones1998efficient}. The resulting expression is
\begin{equation}    \mathrm{EI}_n(\textbf{x}^*)=\left[\Delta_n(\textbf{x}^*)\right]^{+}+\sigma_n(\textbf{x}^*) \varphi\left(\frac{\Delta_n(\textbf{x}^*)}{\sigma_n(\textbf{x}^*)}\right)-\left|\Delta_n(\textbf{x}^*)\right| \Phi\left(\frac{\Delta_n(\textbf{x}^*)}{\sigma_n(\textbf{x}^*)}\right),
\end{equation}
where $\Delta_n(\textbf{x}^*):=f_n^* - \mu_n(\textbf{x}^*)$ is the expected difference in quality between the previous best and the proposed point $\textbf{x}^*$.

The expected improvement algorithm then evaluates at the point with the largest expected improvement,
\begin{equation}
    \textbf{x}_{n+1}=\operatorname{argmax} \mathrm{EI}_n(\textbf{x}^*).
\end{equation}

\subsubsection{BBO realization}
In our problem, evaluation of $f$ is quite expensive since a single run of simulation is time-consuming. So we adopt BBO to enable parallel simulation. In BBO, the acquisition function is slightly different from what we see in \ref{acquisition function}. Here we consider the choice of selecting
$\textbf{x}_{t, k}$, the k-th element of the t-th batch. The maximization-penalization strategy selects $\textbf{x}_{t, k}$ as
\begin{equation}
    \textbf{x}_{t, k}=\arg \max _{\textbf{x} \in \Theta}\left\{G\left(\alpha\left(\textbf{x} ; \mathcal{I}_{t, 0}\right)\right) \prod_{j=1}^{k-1} \varphi\left(\textbf{x} ; \textbf{x}_{t, j}\right)\right\},
    \label{equation:next points}
\end{equation}  
where $\varphi\left(\textbf{x} ; \textbf{x}_{t, j}\right)$ are local local penalizers centered at $\textbf{x}_{t, j}$ and $G: \mathbb{R} \rightarrow \mathbb{R}^{+}$is a differentiable transformation of $\alpha(\textbf{x})$ that keeps it strictly positive without changing the location of its extrema. We will use $G(z)=z$ if $\alpha(\textbf{x})$ is already positive and the soft-plus transformation $G(z)=\ln \left(1+e^z\right)$ elsewhere. This does not require re-estimation of the GP model after each location is selected, just a new optimization of the penalized utility. \cite{gonzalez2016batch}

\noindent
\subsection{LNS-BO Algorithm}
In last subsection, we introduce how BO is adopted to enhance neighborhood search in LNS. In this part we give the pseudo code of LNS-BO algorithm in Algorithm \ref{LNS-BO algorithm}.
As a part of LNS-BO algorithm, the pseudo code of BO enhanced neighborhood search (BOENS) is given in Algorithm \ref{BOENS}.
The pseudo code is based on the LNS framework, and solves the a subproblem in “Repair” with BBO. Here we denote $g(\mathbb{X})$ as the total demand loss of BSS location solution $\mathbb{X}$.

In real situation, the “Repair” process can add multiple BSSs to $\mathbb{X}^t$ each time, while we add one BSS to $\mathbb{X}^t$ in \ref{LNS-BO algorithm}. This consideration of adding one BSS each time is due to the higher sampling efficiency of BO in low dimensional searching space.

\begin{algorithm}[hbt!]

\caption{BOENS}\label{BOENS}
\KwIn{Number of evaluation for each iteration in BBO $m$, times of sampling $n\_{sample}$, solution $x_0$, evaluated points $\mathbf{X}_{step}$ with objective function value $\mathbf{Y}_{step}$, unrepaired BSS location solution $\mathbb{X}^t$}

\KwOut{Location record $\mathbf{X}_{step}$, objective function value record $\mathbf{Y}_{step}$}
$current\_iter \leftarrow 0$ \;
\While{$current\_iter < n\_{sample}$} {
            Construct the GP model according to (\ref{equation:update bo1}) and (\ref{equation:update bo2}) with the observations $\mathbf{X}_{step}$ and $\mathbf{Y}_{step}$ \;
            Suggest the next $m$ points to sample by (\ref{equation:next points}), and evaluate them\;
            Let $\mathbf{X}_{s}$ be the set of sampled points, and $\mathbf{Y}_{s}$ be the set of their objective values\;
            Update  
            $\mathbf{X}_{step} \leftarrow \mathbf{X}_{step} \cup \mathbf{X}_{s}$, $\mathbf{Y}_{step} \leftarrow \mathbf{Y}_{step} \cup \mathbf{Y}_{s}$   \;
        
            
           
            $current\_iter \leftarrow current\_iter + 1$  \;
            
        }

\end{algorithm}

\begin{algorithm}[hbt!]
\caption{LNS-BO} \label{LNS-BO algorithm}
\KwIn{Number of BSS destruction $k$, number of initial evaluation in BBO $n$, number of evaluation for each iteration in BBO $m$, times of sampling $n\_{sample}$, initial BSS location solution $\mathbb{X}^0$}

\KwOut{Final BSS location solution $\mathbb{X}^{b}$}

$\mathbb{X}^b \leftarrow \mathbb{X}^0$, $\mathbb{X} \leftarrow \mathbb{X}^0$ \;

$last\_obj \leftarrow g(\mathbb{X}^0)$, $best\_obj \leftarrow g(\mathbb{X}^0)$ \;

\While{not terminate} {
    \tcc{Destroy}
    Randomly select $k$ stations in $\mathbb{X}$ and remove them, obtaining $\mathbb{X}^{t}$  \;

    \tcc{Repair}
    \For{$i \leftarrow 1$ \KwTo $k$  }  {
        Randomly sample $n$ points $\mathbf{X}_{init} = \left\{ \textbf{x}_1, \textbf{x}_2, \ldots, \textbf{x}_n \right\} $ for the position of the $i$-th BBS \;
        Let $ \mathbb{X}^{tmp}_j = \left\{ \textbf{x}_j \right\} \cup \mathbb{X}^t , j \in \{1,2,\ldots,n\}$  and evaluate $\mathbb{X}^{tmp}$ to obtain $\mathbf{Y}_{init} = \Big\{ g( \mathbb{X}^{tmp}_1 ),g( \mathbb{X}^{tmp}_2 ),\ldots,g( \mathbb{X}^{tmp}_n ) \Big\}$     \tcp*[r]{Parallel computing}

        $\mathbf{X}_{step} \leftarrow \mathbf{X}_{init}$ \;
        $\mathbf{Y}_{step} \leftarrow \mathbf{Y}_{init}$ \;
        
        $\mathbf{X}_{step}, \mathbf{Y}_{step} \leftarrow$ BOENS{$\Big( \mathbb{X}^t, \mathbf{X}_{step},\mathbf{Y}_{step},n\_{sample},m \Big) $} \;

        $j^* \leftarrow \arg \min_{j}  \mathbf{Y}_{step,j} $ \;

        $obj \leftarrow \mathbf{Y}_{step,j^*}$ \;
        $\mathbb{X}^t \leftarrow \{ \mathbf{X}_{step,j^*} \} \cup \mathbb{X}^t$   \;             
        }

        \tcc{Update solutions}
        \If{$obj<last\_obj$}   {
            $\mathbb{X} \leftarrow \mathbb{X}^t$  \;  
            $last\_obj \leftarrow obj$  \;
        }

        \If{$obj < best\_obj$}{
            $best\_obj \leftarrow obj$ \;
            $\mathbb{X}^b \leftarrow \mathbb{X}^t$ \;
        }

    }
    
\end{algorithm}

\section{Numerical experiment}
\label{section:Num_Results}

In this section, we first build a simulation model in accordance with MILP model, validate their consistency and find out the disadvantage of mathematical model in \ref{exp:1}. Then we make a comparison between our simulation optimization approach and traditional mathematical method in solving different scale location problem in \ref{exp:2}, where the effectiveness of LNS-BO over SA is also tested. In \ref{exp:3}, we apply simulation optimization method to a real instance, and demonstrate its advantages in solving large scale problems.

\subsection{Disadvantage of mathematical model} \label{exp:1}

\begin{figure} [hbt!]
    \centering
    \includegraphics[width=0.9\linewidth]{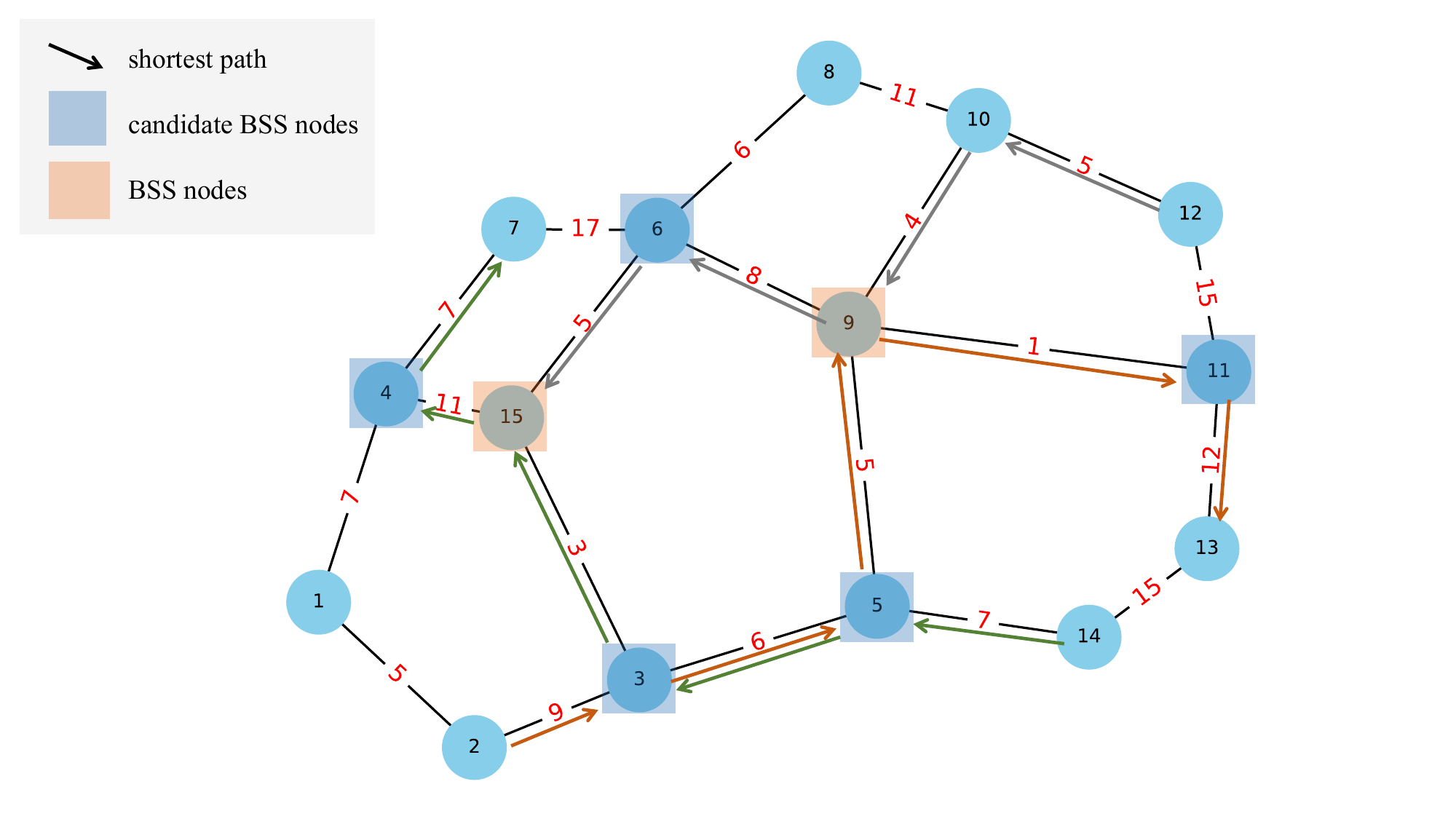}
    \caption{Testing map for part \ref{exp:1}}
    \label{fig:math-case-map}
\end{figure}

In this section, we explore to what extent will inaccurate AWT influence optimal BSS location selection(solved  by MILP model) compared to actual optimal BSS location selection (solved by simulation optimization). 

We set up a testing map as shown in Figure\ref{fig:math-case-map}. There are 15 dots, 3 path demands, 2 BSSs and 7 candidate BSS nodes. The number on the black line indicates the traveling cost between two dots. Three path demands are \{(2,13), (12,15), (14,7)\}, whose shorted paths are illustrated by consecutive arrows on the map. We use 7 squares, say \{3, 4, 15, 6, 9, 11, 5\}, to represent 7 candidate BSS nodes, and two orange ones stand for nodes where BSSs are built. Possible initial battery levels of the vehicles are \{30, 35, 40, 45, 50, 55, 60, 65, 70, 75, 80, 85, 90, 95, 100\}.

For parameters, we set $\lambda_{(2,13)}=3600/210$, $\lambda_{(12,15)}=3600/280$, $\lambda_{(14,7)}=3600/300$. We set $d_{max}=5$, $c=2$, $\epsilon=0.3$, $N=2$, $\mu_j=3600/300$ and $LL=7$. To compute $\pi^b_{j,(i_s,i_e)}$ (see \ref{pi}), we set $\alpha_1=0.04$, $\alpha_2=0.8$, $\alpha_3=0.02/60$. In MILP model, we set $l_3=345/3600$, $l_4=312/3600$, $l_{15}=415/3600$, $l_6=380/3600$, $l_9=612/3600$, $l_11=200/3600$, $l_5=250/3600$, where $l_j$ is estimated artificially, indicating a flaw of MILP model. However, in simulation optimization model, $l_j$ is computed according to historical queuing information.

In the experiment, fistly we get the optimal BSS location  solution [11, 5] by solving MILP model under our experiment setup, with $loss=44.5\%$ in accurate estimator. Using the same parameters except AWT, we obtain the optimal BSS location solution [11, 15] by enumeration, with $loss=40.8\%$ in accurate estimator. So the true optimal solution has a lower $loss$ of 3.7\% compared to MILP solution. Thus we conclude that inaccurate AWT will cause inaccurate objective function value evaluation, leading to suboptimal BSS location selection.

\subsection{Comparison between simulation optimization and MILP}  \label{exp:2}
In this section, we will discuss the performance of MILP model and simulation optimization model under different road network scales and different demand sizes, and evaluate the optimization effect of LNS-BO algorithm compared with SA algorithm. The experiment reasonably cuts the Shanghai road network to obtain three different scales of road networks with 50, 101 and 154 map nodes. Road networks of different scales can better simulate the actual road network conditions of various complexities; the experiment also considers the impact of different demands. There are three different scenarios of low demand, medium demand and high demand for road networks of different scales, which avoids the contingency caused by the combination of special demand numbers and special demand locations, and comprehensively evaluates the performance of the algorithm at different demand scales.

These experiments set the same simulation running time, which is two hours, to compare the optimization effects of different models and algorithms horizontally. The final experimental data results are shown in Table\ref{tab:exp2}. 

The data shows that the mathematical model can obtain a solution within a limited time when dealing with small road networks and small-scale demand numbers. For large road networks and large-scale demand numbers, the mathematical model cannot even obtain a feasible solution within a limited time; and because the mathematical model makes random assumptions about the average waiting time for each BSS, the optimal or better value obtained under some conditions is even not as good as the initial value of the simulation optimization. As for simulation optimization, SA and LNS-BO can handle scenarios with any road networks and path demand numbers. They give better solutions in a limited time, demonstrating the advantages of simulation models. However, we notice that LNS-BO is generally just sightly better than the optimization results of the SA algorithm. We speculate that these instances are not complex enough, even random neighborhood search could be an effective approach. When faced with real world instance, the advantages of LNS-BO will be shown, which we illustrate in the next experiment.

\begin{table}[hbt!] 
\caption{Comparison between MILP and simulation optimization on different instances.}
\resizebox{\textwidth}{!}{
\begin{threeparttable}[b]
\begin{tabular}{cccccc}
\toprule[1.5pt]
\multirow{2}{*}{\textbf{Traffic Network Size}} & \multirow{2}{*}{\textbf{Path Demand Number}} & \multirow{2}{*}{\textit{MILP}} & \multicolumn{3}{c}{Simulation Optimization}       \\ \cline{4-6} 
                                               &                                              &                                & Initial Obj\tnote{1} & \textit{SA}      & \textit{LNS-BO}    \\ \midrule[1.1pt]
\multirow{3}{*}{50}                            & small                                        &  {1.13\%}                & 5.94\%      &  {0\%}  &  {0.01\%}     \\
                                               & medium                                       &  {8.79\%}                & 16.20\%     &  {5.31\%}  &  {4.76\%}  \\
                                               & large                                        &  {47.91\%}               & 47.80\%     &  {44.45\%} &  {44.20\%} \\ \hline
\multirow{3}{*}{101}                           & small                                        &  {2.76\%}                & 9.15\%      &  {0.87\%}  &  {0.08\%}  \\
                                               & medium                                       &  {8.45\%}                & 14.35\%     &  {2.94\%}  &  {2.44\%}  \\
                                               & large                                        &  {/} \tnote{2}                    & 49.29\%     &  {42.04\%} &  {43.67\%} \\ \hline
\multirow{3}{*}{154}                           & small                                        &  {11.68\%}               & 3.43\%      &  {1.11\%}  &  {0.73\%}  \\
                                               & medium                                       &  {/}                     & 23.11\%     &  {14.57\%} &  {14.56\%} \\
                                               & large                                        &  {/}                     & 48.20\%     &  {45.97\%} &  {45.91\%} \\ \bottomrule[1.5pt]
\end{tabular}
\begin{tablenotes}
     \item[1] SA and LNS-BO are used in Simulation Optimization method, where Initial Obj stands for the initialized objective value before optimization.
     \item[2] ``/" stands for ``No feasible solution found in given time".
   \end{tablenotes}
\end{threeparttable}
}
\label{tab:exp2}
\end{table}



\subsection{simulation optimization applied to Large scale problem} \label{exp:3}

The dataset utilized in our research is the 2007 Shanghai GPS taxi dataset, specifically from July 20th. The dataset includes fields such as taxi ID, timestamp, latitude, longitude, speed (km/h), heading angle, and passenger status. There is a total of 6,075,588 records in this dataset, with an average sampling frequency of 60.5 seconds. Before we carry out this experiment, we need to preprocess the data, extract path demand, extract and build road network and configure simulation model. These are illustrated in appendix \ref{appendix:exp3}. Considering the feasibility of station placement in the real world, candidate BSS nodes include those where supermarkets, parking lots, charging stations, malls, or BSSs have already been established. The simulation length is set as 10 days. By preliminary tests, we have validated that the simulation length is long enough to guarantee a small output variation. More specifically, on the current NEO BSS layout, we have performed 10 replications of the simulation. The averaged obejective value is 63.15\%, and the standard deviation is 0.12\%.

\begin{figure}[hbt!]
    \centering
    \includegraphics[width=0.7\linewidth]{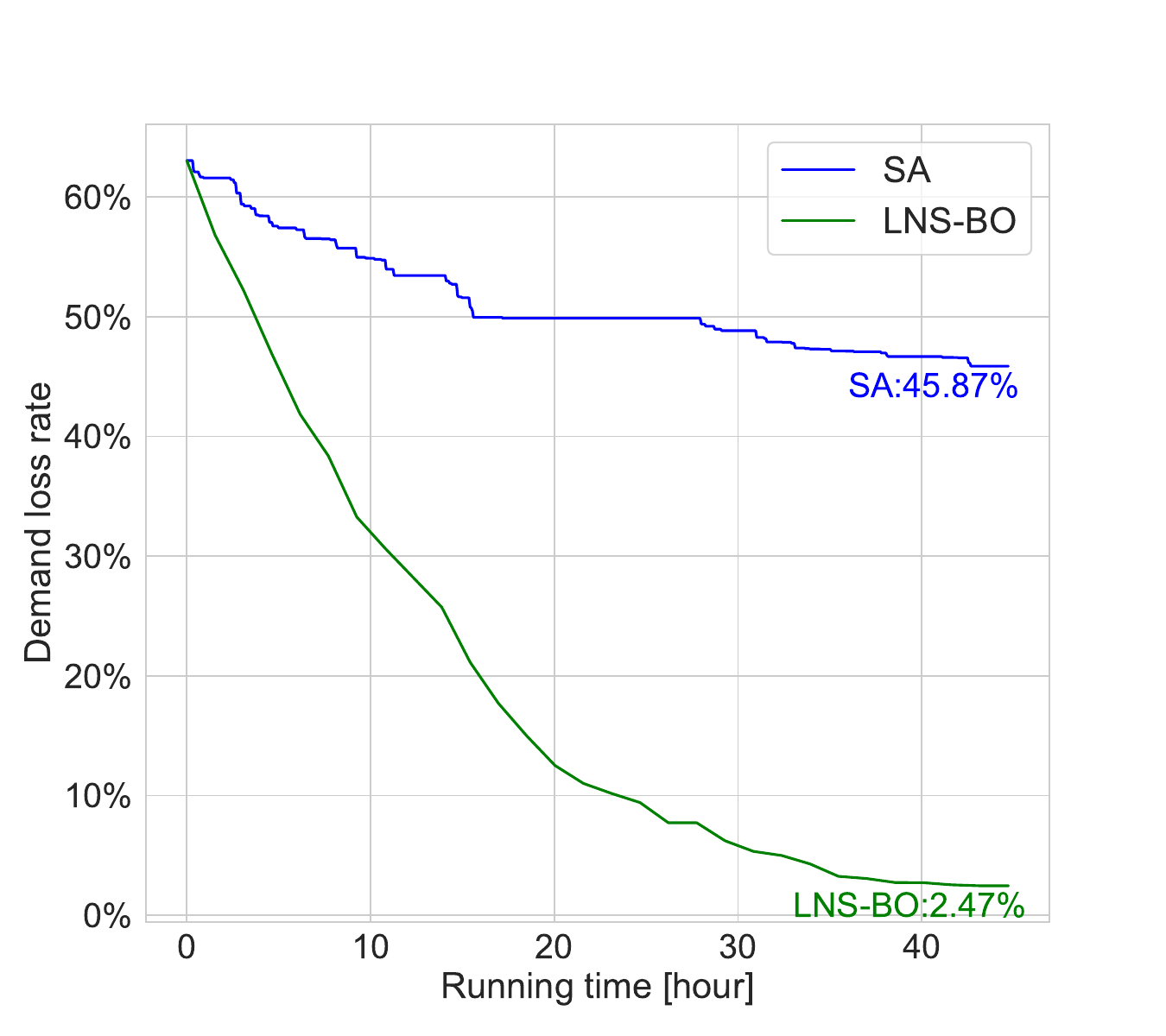}
    \caption{Comparison of SA optimized BSSs layout and LNS-BO optimized BSSs layout. The initial layout (at Running Time 0) is current NIO BSS layout in Shanghai, China}
    \label{fig:new old-comparison}
\end{figure}


Then we make a comparison between the SA and LNS-BO. With the same computation time, the optimization results for SA and LNS-BO algorithms are shown in the Figure\ref{fig:new old-comparison}. By the 44th hour, the demand loss rate of LNS-BO decreases to 2.47\%, while at this point, SA's optimization result is just as good as what LNS-BO achieved in the 5th hour. Furthermore, up to the 44th hour, SA still maintains an 45.87\% demand loss rate. This demonstrates that LNS-BO has significantly improved the optimization efficiency and effectiveness compared to SA, especially in high-cost optimization problems.

Additionally, the initial layout of both SA and LNS-BO uses the current NIO BSS layout in Shanghai, China, whose demand loss rate is 75.6\% (as indicated by the starting point in figure\ref{fig:new old-comparison}). This suggests that there are still certain issues with the layout of NIO BSSs in Shanghai under the current experimental data, leaving substantial room for further optimization.

\begin{figure}[hbt!]
    \centering
    \subfigure[\centering Before optimization. 
     Demand loss rate: 63.04\%]{\includegraphics[width=0.62\linewidth]{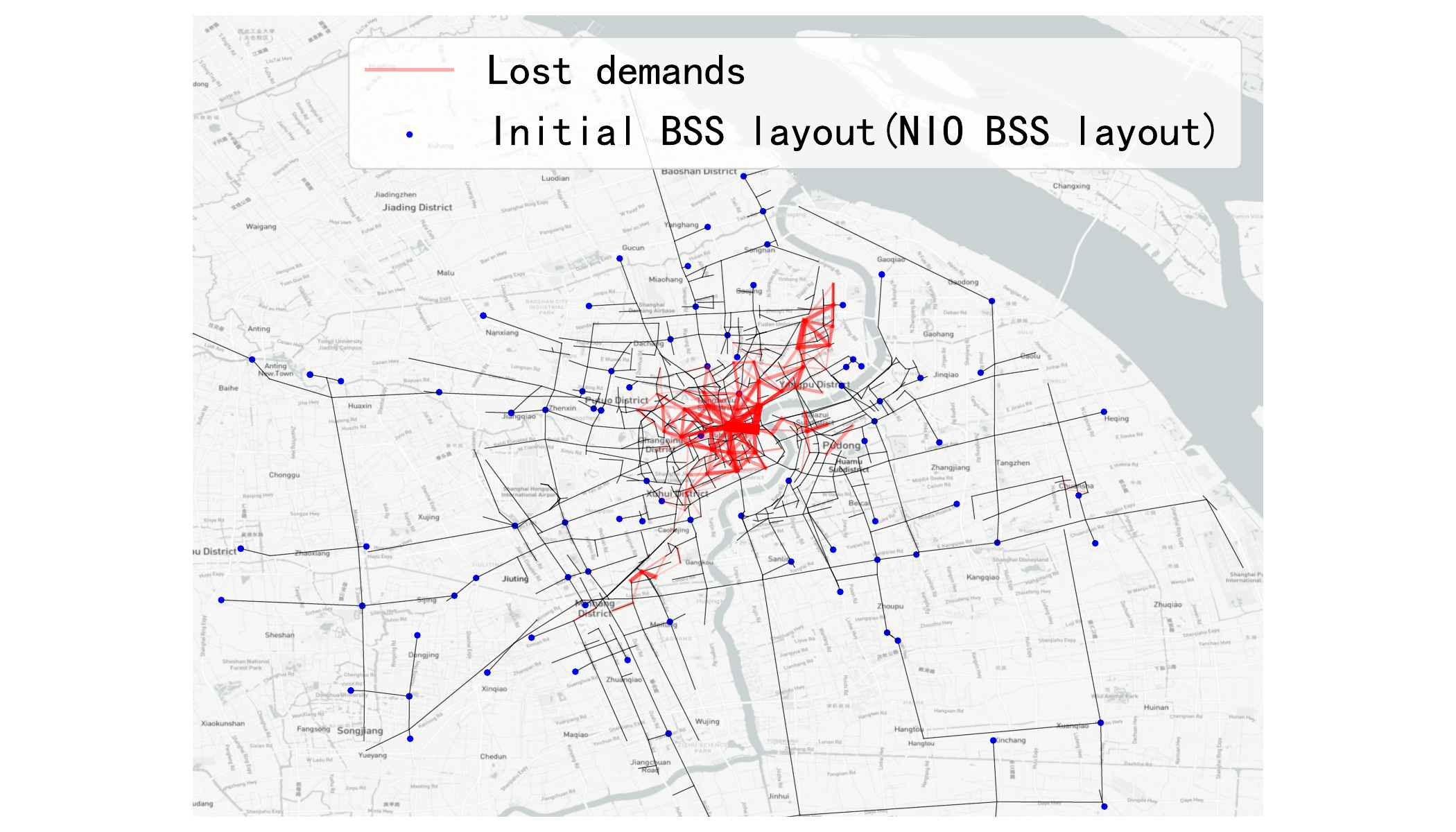}}
    \subfigure[After SA optimization. 
      Demand loss rate: 45.87\%]{\includegraphics[width=0.64\linewidth]{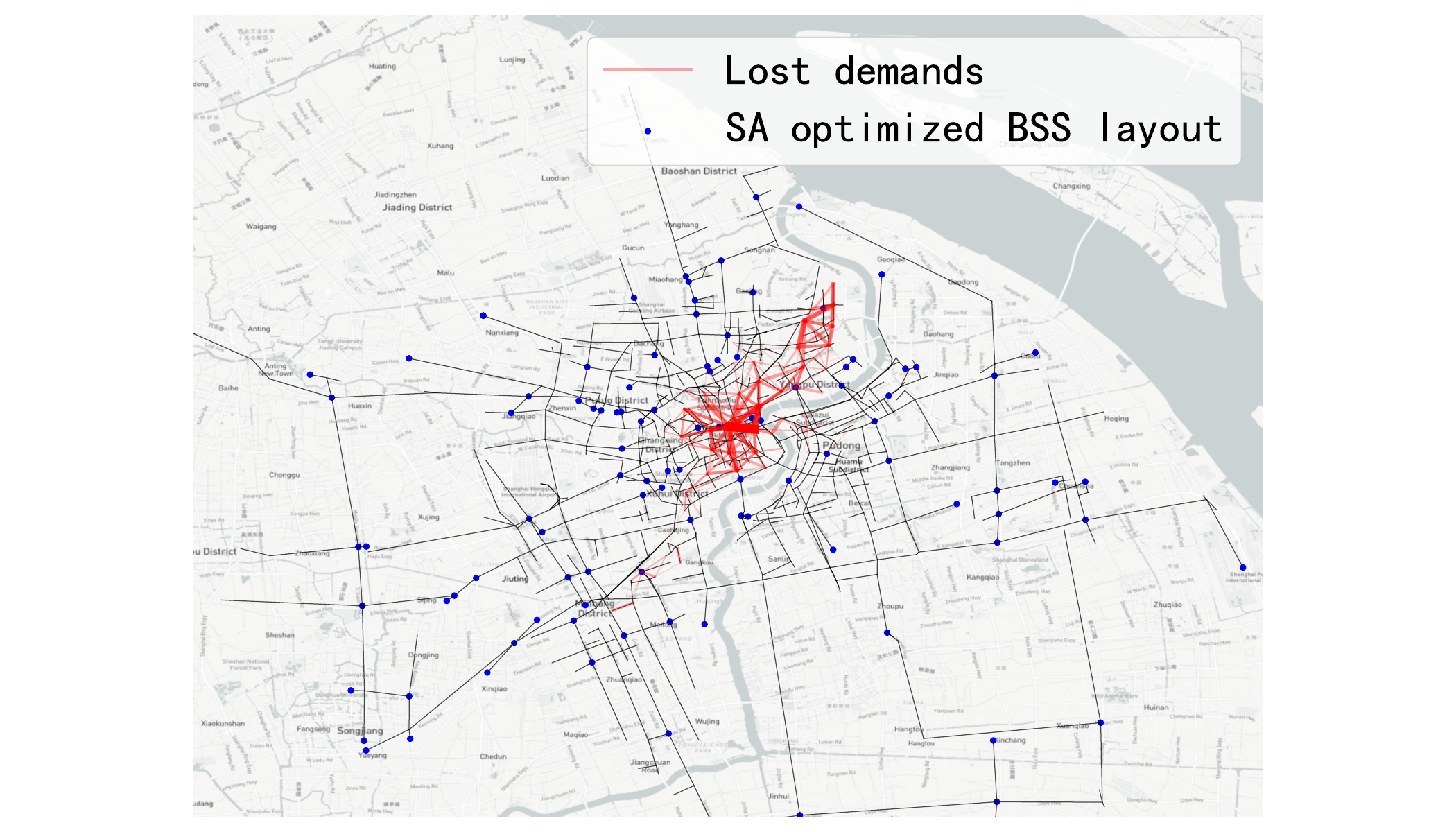}}
    \subfigure[After LNS-BO optimization.
     Demand loss rate: 2.47\%]{\includegraphics[width=0.64\linewidth]{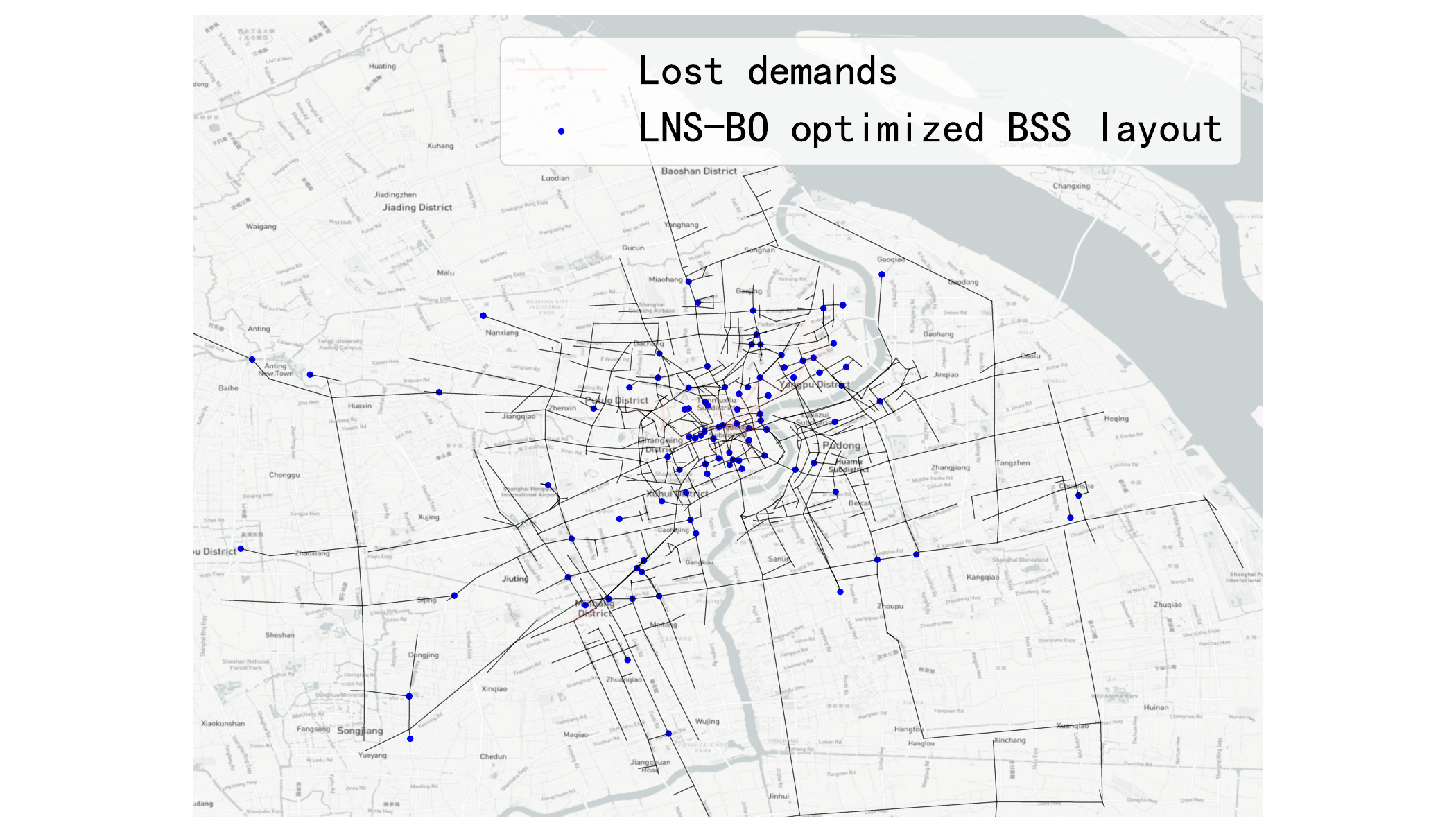}}
    
    \caption{Comparison of initial BSSs layout (NIO BSS layout), SA optimized BSS layout and LNS-BO optimized BSS layout}
    \label{fig:layout comparison}
\end{figure}

From Figure\ref{fig:layout comparison}, visually, NIO's BSSs layout appears relatively uniform. However, we can see large demand loss in central area of Shanghai, meaning BSSs placed in central area are not sufficient. After SA optimization, it's apparent that BSSs are driven towards central districts, decreasing demand loss rate pretty decently. In the case of LNS-BO, there're essentially no lost demands after optimization. We can see that many BSSs concentrate in the central area, absorbing huge demands there.

\section{Conclusion}
\label{section:conclusion}

The location selection problem of EV BSSs is of great significance. This paper first starts from the mathematical model and establishes the MILP model of this problem based on the latest research. Since the average waiting time of the BSS in the MILP model is a hyperparameter, there are certain problems in the evaluation of its objective function. We compare it with the simulation model on a small case to verify this conclusion, and this defect will cause the location selection plan finally solved by the mathematical model to be suboptimal. Furthermore, in terms of performance, we found that as the scale of the problem increases, even with the most advanced solver, it is still difficult for the mathematical model to find a feasible solution within the specified time, but the simulation optimization method solves this problem well. Besides, under the framework of simulation optimization, we proposed an efficient algorithm, which is far better than the traditional simulated annealing algorithm in both optimization speed and optimization effect. This has been verified not only in the manually set experiment, but also in the actual city case of Shanghai, which well verifies the efficiency of the new algorithm.

\section*{Acknowledgements}

This work was supported by the Student Innovation Training Program of Tongji University under Grant [X2022191]; Tongji University Fundamental Research Funds for the Central Universities under Grant [22120240607].

\appendix

\section{Procedures of experiment 5.3} \label{appendix:exp3}
To introduce what the dataset looks like, part of the data is shown in Table \ref{data table}.

\begin{table}[htb!]     
    \centering
    \resizebox{0.9\textwidth}{30mm}{
    \begin{tabular}{rccccccc}
    \toprule[1.5pt]
    \multicolumn{1}{c}{\textbf{}} & \textbf{VehicleNum} & \textbf{Time} & \textbf{Longtitude} & \textbf{Latitude} & \textbf{Speed} & \textbf{Degrees} & \textbf{Status} \\ \midrule[1pt]
    0                & 10001 & 2007-02-20 00:02:27 & 121.423167 & 31.165233 & 7.0 & 116.0 & 3.0 \\
    \textit{1}       & 10001 & 2007-02-20 00:05:36 & 121.423167 & 31.165233 & 7.0 & 116.0 & 3.0 \\
    \textit{2}       & 10001 & 2007-02-20 00:08:45 & 121.423167 & 31.165233 & 7.0 & 116.0 & 3.0 \\
    \textit{3}       & 10001 & 2007-02-20 00:11:55 & 121.423167 & 31.165233 & 7.0 & 116.0 & 3.0 \\
    \textit{4}       & 10001 & 2007-02-20 00:15:04 & 121.423167 & 31.165233 & 7.0 & 116.0 & 3.0 \\
    ...              & ...   & ...                 & ...        & ...       & ... & ...   & ... \\
    \textit{6075583} & 99780 & 2007-02-20 23:54:10 & 121.480600 & 31.241300 & 0.0 & 45.0  & 0.0 \\
    \textit{6075584} & 99780 & 2007-02-20 23:56:13 & 121.480800 & 31.241300 & 0.0 & 45.0  & 0.0 \\
    \textit{6075585} & 99780 & 2007-02-20 23:58:15 & 121.480600 & 31.241300 & 0.0 & 45.0  & 0.0 \\
    \textit{6075586} & 99780 & 2007-02-20 23:59:15 & 121.480600 & 31.241300 & 0.0 & 45.0  & 0.0 \\
    \textit{6075587} & Nan   & NaN                 & NaN        & NaN       & NaN & NaN   & NaN \\ 
    \bottomrule[1.5pt]
    \end{tabular}}
    \caption{GPS data of taxis, Shanghai}
    \label{data table}
\end{table}

\subsection{Data preprocessing}
There are some data in the original data that interfere with our research, and data preprocessing must be performed before the experiment begins.
\begin{enumerate}
    \item Remove the tuples containing null values, get 6,075,587 pieces of data, and filter out 1 piece of data.
    \item Limit the spatial range of GPS data points to Shanghai. The longitude range is [120.852326, 122.118227] and the latitude range is [30.691701, 31.874634]. 6,068,426 pieces of data were obtained, and 7,161 pieces of data were filtered out.
    \item Filter out the data with instantaneous changes in passenger status, 6,021,469 pieces of data are obtained, and 46,957 pieces of data are filtered out.
    \item Then extract a total of 745,832 pieces of data in the passenger-carrying state, and extract the data paths in the passenger-carrying state, and obtain a total of 66,431 paths to prepare for the subsequent extraction of path requirements and road network extraction.
\end{enumerate}

\subsection{Path demand extraction}
After Gridding the previously obtained paths, we obtained 20,124 path requirements. 

\begin{figure}
    \centering

    \begin{tikzpicture}
        \node[anchor=south west] at (0, 3.4) {\includegraphics[width=4.7cm,height=3.4cm]{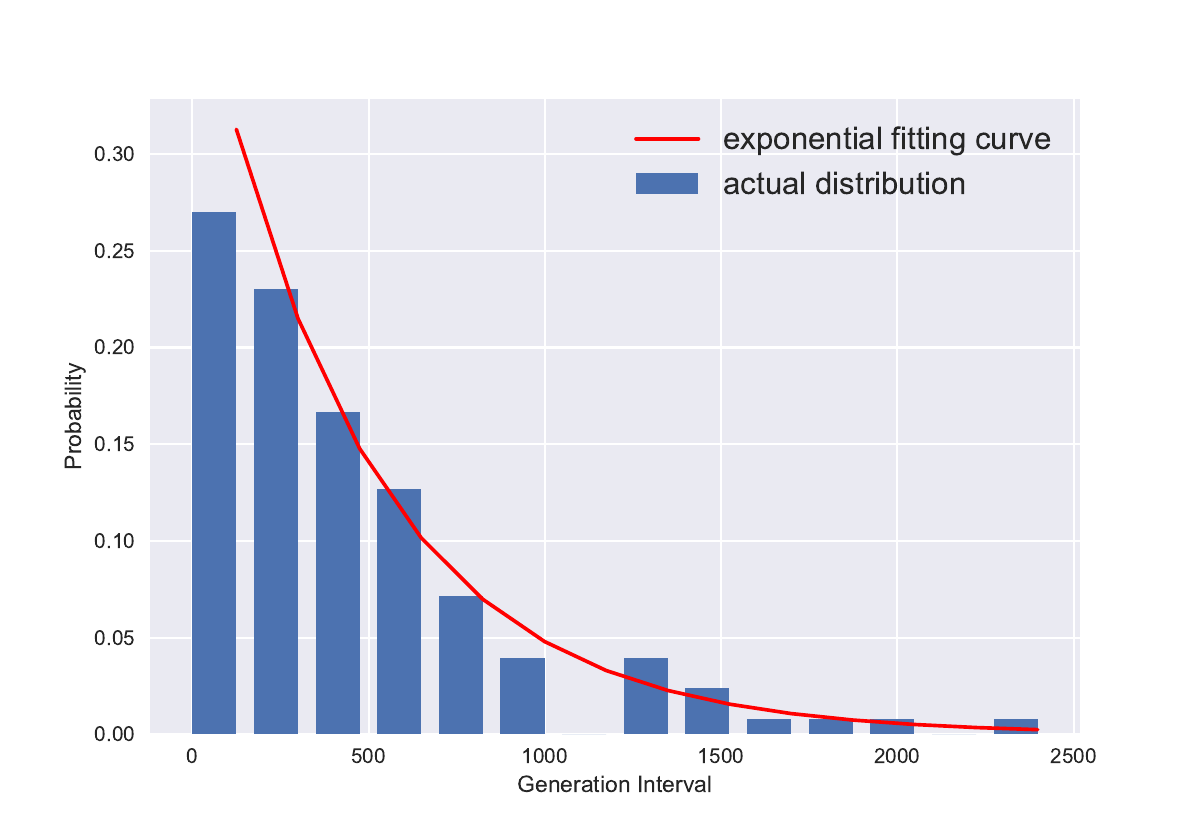}};
        \node[anchor=south west] at (0, 0.23) {\includegraphics[width=4.7cm,height=3.4cm]{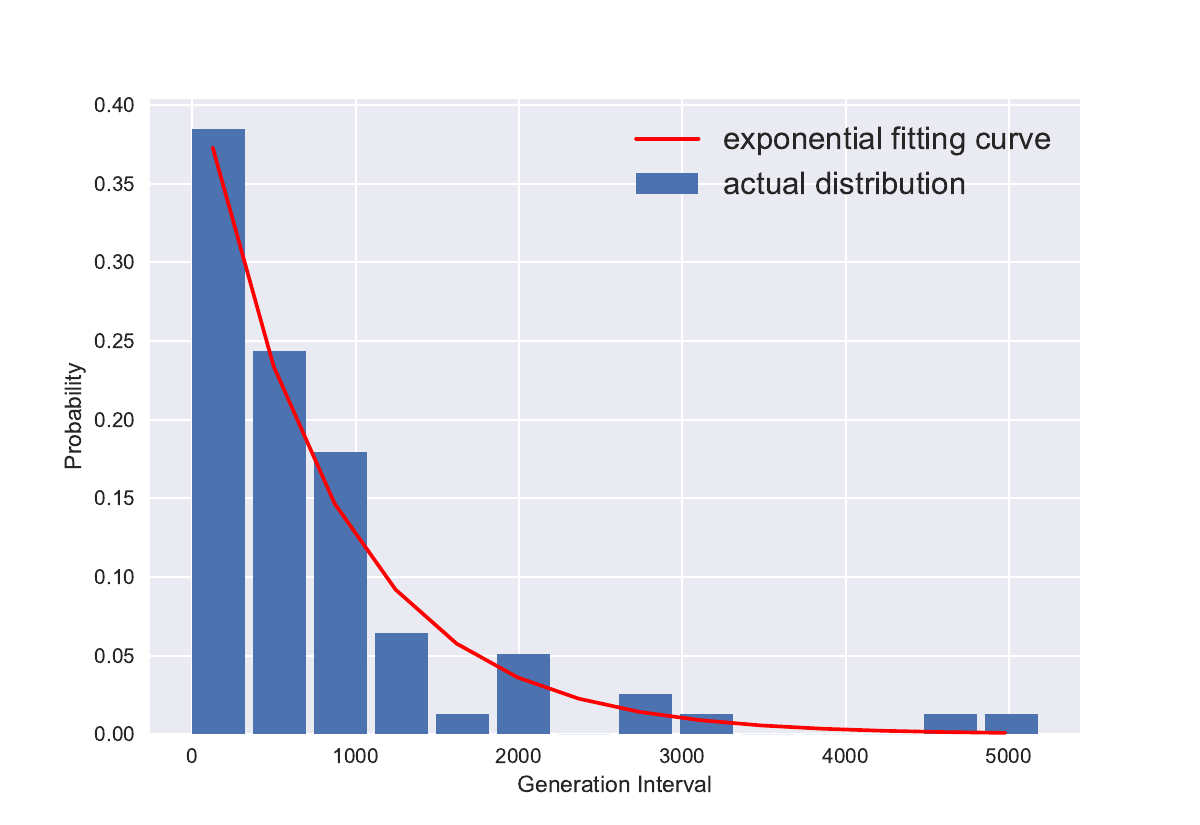}};
        
        \node[anchor=south west] at (4.5, 0) {\includegraphics[width=10cm,height=7cm]{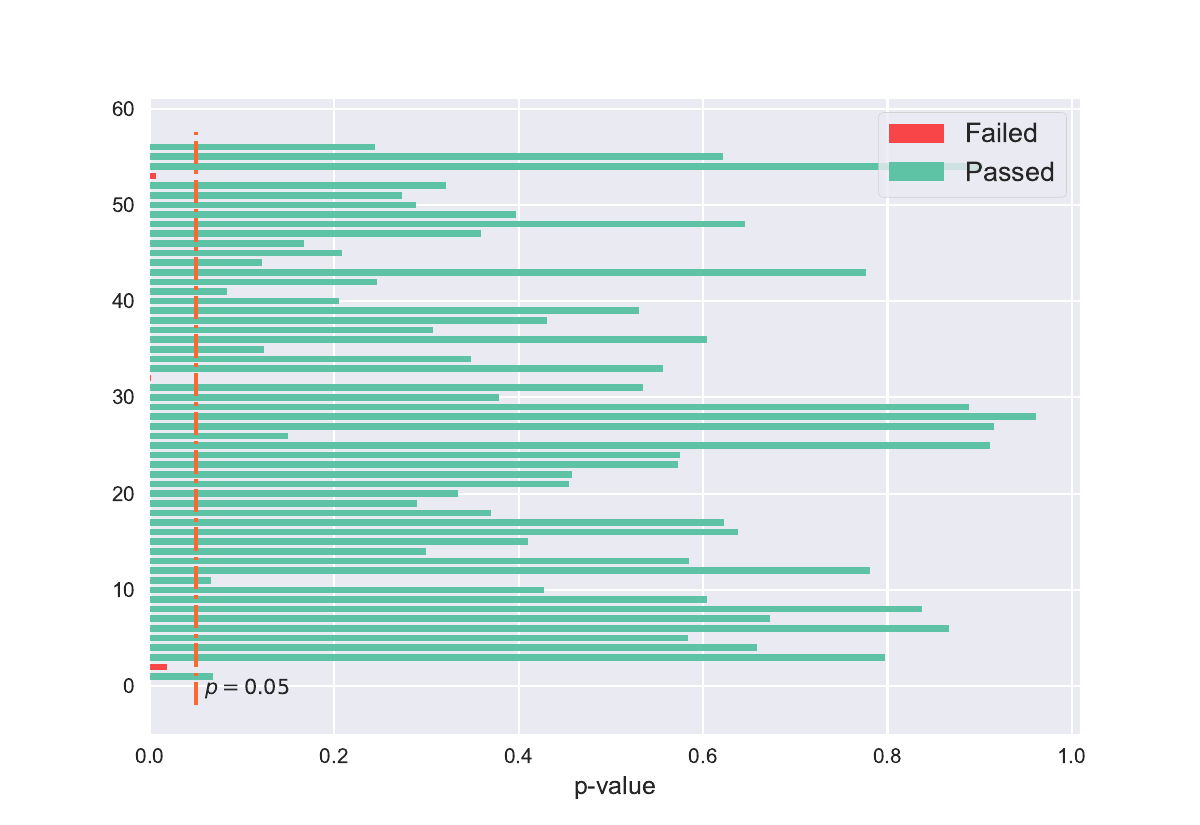}};
        
        \node at (3, 0) {\fontsize{6}{6}\selectfont Path demand fitting}; 
        \node at (9.7, 0) {\fontsize{6}{6}\selectfont Hypothesis testing};  
    \end{tikzpicture}

    \caption{Path demand generation distribution modeling}
    \label{fig:hypothesis testing}
\end{figure}

In order to model the distribution of path demand generation time interval, we extracted a total of 56 path demand with demand greater than 50 and plotted their frequency-generation interval histogram as shown in the left graphs of figure\ref{fig:hypothesis testing}(2 random path demands). From the graphs, We discovered that exponential distribution may fit the generation interval distribution well. Thus, we did hypothesis testing on our assumption. As shown in right graph of figure\ref{fig:hypothesis testing}, the p value corresponding to 53 path demand was greater than 0.05, and the pass rate was 94.6\%, so we accepted the null hypothesis, indicating that the exponential distribution has a good fitting effect. Therefore, we used the exponential distribution to fit the path demand generation interval, and obtained the parameters of the exponential distribution probability density function for use in the flexsim simulation software.




\subsection{Road network extraction}
In order to further use the real-life map in the Flexsim simulation software, it is necessary to extract the nodes and edges of Shanghai City. The method that we use is shown in Figure\ref{fig:路网生成步骤}. 
The final transportation network diagram is shown in figure\ref{fig:k_means_1000_map} below. The red dots represent key points, i.e. crossroads and end points. The blue line stands for the road on which the number represents the length of the road segment in reality. Notice that here the transportation network is a abstract concept, since we only illustrate the key points and their connection. And we incorporate the trajectory characteristics in the real lengths of road segments, which is what we need in our problem.

\begin{figure}[hbt!]
    \centering
    \includegraphics[width=\linewidth]{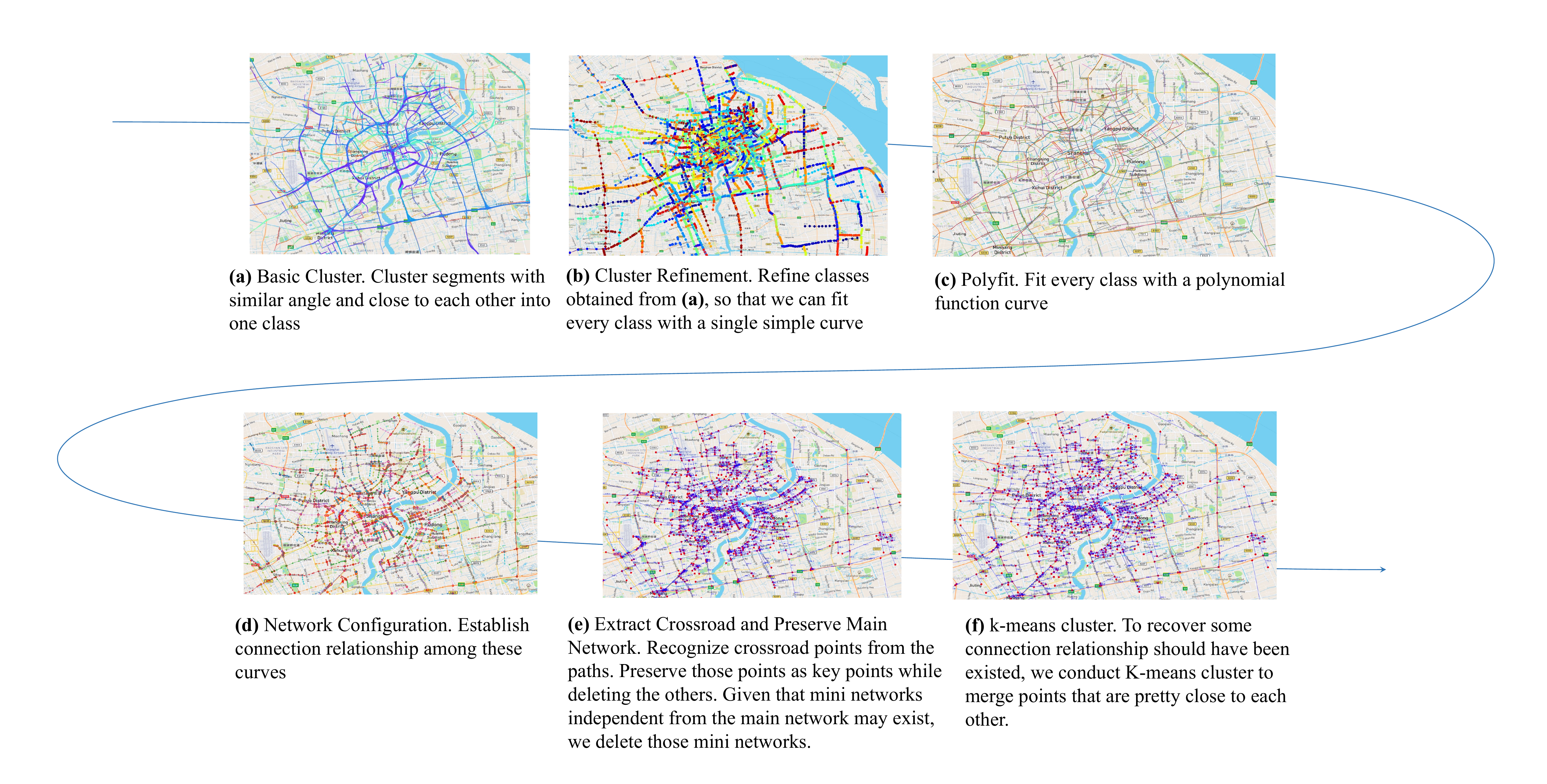}
    \caption{Road network extraction through GPS data}
    \label{fig:路网生成步骤}
\end{figure}

\begin{figure}[hbt!]
    \centering
    \includegraphics[width=0.8\linewidth]{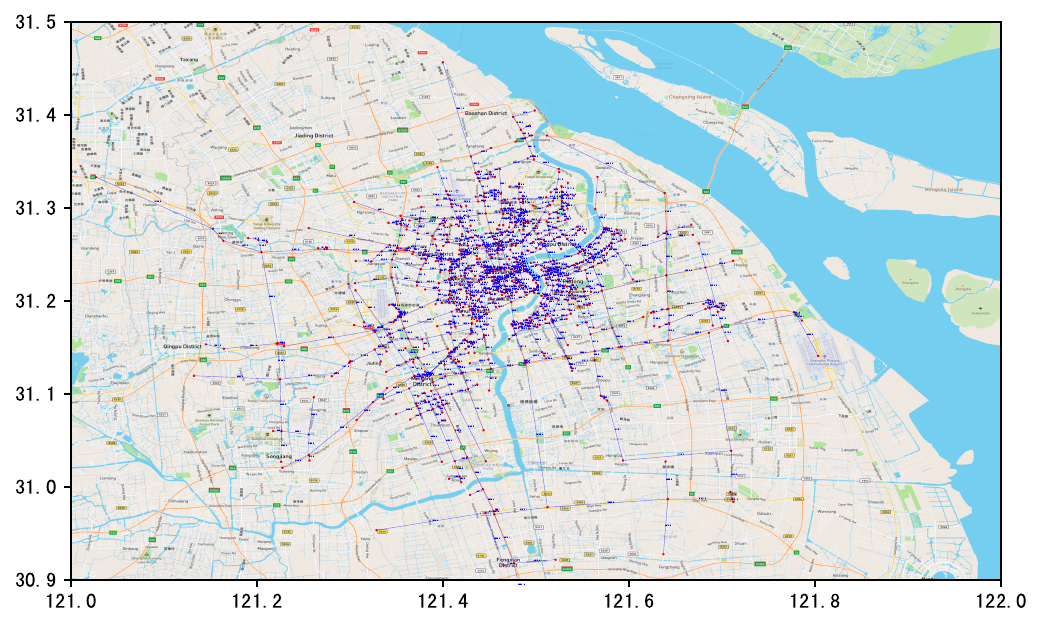}
    \caption{Road network of Shanghai, China}
    \label{fig:k_means_1000_map}
\end{figure}

\subsection{Simulation model generation} \label{simulation}
We import the data extracted from the existing data set as path demand parameters into the simulation model. For each path demand, the distribution of its demand generation follows an exponential distribution with the parameter corresponding to the demand $\lambda$.

For the road networks of different models, we cut them from the actual road network in Shanghai. The distance between the two points abstracted is the length of the actual road network, not the Euclidean distance between the two points.

For the simulation process of each model, the warm-up time is 300,000 seconds, and the simulation time of each model is 864,000 seconds, that is, 10 days. The final objective function result can be obtained after the simulation.

For each battery swap station, the maximum capacity is set to 6. Failures and repairs are not considered in the simulation. The battery swap time of each battery swap station follows an exponential distribution with a mean of 300 seconds.

For all vehicles in the road network, we set the maximum detour distance $S_{max}$ acceptable to users for path demand to 5km, the maximum tolerable battery swap queue length to 6, the initial power $L_c$ of the vehicle follows a uniform distribution of (5,100)Kwh, and the average speed $v_0$ of the vehicle is 10m/s. The battery capacity of the vehicle is 100Kwh, the power consumption per kilometer is 2Kwh, and the final battery swap station selection process for each vehicle conforms to the formula mentioned in section\ref{problem description}.

\bibliography{mybibfile}

\end{document}